\documentclass[15pt,fleqn,leqno]{article}
\usepackage{graphicx,color,enumerate,amssymb}
\usepackage[authoryear,round]{natbib}
\usepackage{amsmath}
\usepackage{amssymb}
\makeatletter
\@addtoreset{equation}{section}
\renewcommand\thetable{\@arabic\c@table}
\renewcommand\thefigure{\@arabic\c@figure}
\long\def\@makecaption#1#2{%
  \vskip\abovecaptionskip
  \begin{center}%
  \sbox\@tempboxa{#1: #2}%
  \ifdim \wd\@tempboxa >\hsize
    #1: #2\par
  \else
    \global \@minipagefalse
    \hb@xt@\hsize{\hfil\box\@tempboxa\hfil}%
  \fi
  \end{center}%
  \vskip\belowcaptionskip}
\def\N{{\rm I\kern-.15em N}}
\def\R{{\rm I\kern-.2em R}}
\def\Z{{\rm Z\kern-.26em Z}}
\makeatother

\newtheorem{thm}{Theorem}[section]

\newtheorem{prop}[thm]{Proposition}
\newtheorem{rem}[thm]{Remark}

\newcommand{\be}{\begin{eqnarray}}
\newcommand{\ee}{\end{eqnarray}}
\newcommand{\bq}{\begin{eqnarray*}}
\newcommand{\eq}{\end{eqnarray*}}

\newcommand{\bewend}{\hspace*{2mm}\rule{3mm}{3mm}}

\begin{document}
\begin{center}
{ \Large MLE--EQUIVARIANCE, DATA TRANSFORMATIONS AND  INVARIANT TESTS OF FIT  \footnote{Proposed running head:
MLE--equivariance, data transformations, and invariant tests}
} \vspace*{0.5cm} {\large\sc }\\ \vspace*{0.5cm}  {\large\sc Muneya Matsui$^{(1)}$\footnote{
Muneya Matsui's research is partly supported by the JSPS Grant-in-Aid for Scientific Research C (23K11019).}, Simos Meintanis$^{(2),(3)}$} \\
{\it (1) Department of Business Administration, Nanzan University,
Nagoya, Japan},\footnote{Email-address: mmuneya@gmail.com} \\
{\it (2) Department of Economics, National and Kapodistrian University
of Athens, Athens, Greece}\footnote{On sabbatical leave from the University of Athens, Email-address: simosmei@econ.uoa.gr} \\
{\it (3) Pure and Applied Analytics, North--West University
, Potchefstroom, South Africa} \\
\end{center}

\vspace*{1cm} {\small {\bf Abstract.} We define  data transformations that leave certain classes of distributions invariant, while acting in a specific manner upon the parameters of the said distributions. It is shown that under such transformations the maximum likelihood estimators behave in exactly the same way as the parameters being estimated. As a consequence goodness--of--fit tests based on standardized data obtained through the inverse of this invariant data--transformation reduce to the case of testing a standard member of the family with fixed parameter values. While presenting our results, we also provide a selective  review of  the subject of equivariant estimators always in connection to invariant goodness--of--fit tests. A  small Monte Carlo study is presented for the special case of testing for the Weibull distribution, along with real--data illustrations.}\\

\vspace*{0.3cm} {\small {\it Keywords.} MLE--equivariance;   Goodness--of--fit test; Weibull distribution; Test for exponentiality
\\ } 
\newpage

\section{Introduction}\label{sec_1}

\vspace*{-.3cm} Suppose that $\widehat \theta_n:=\widehat\theta_n(X_1,...,X_n)$ is the maximum likelihood estimator (MLE) of a parameter (vector) 
$\theta$ based on independent copies $(X_j, \ j=1,...,n)$ of a random variable $X$ with distribution $F$.  Then by the invariance property of the MLE, 
$h(\widehat \theta_n)$ is the MLE of the parameter $h(\theta)$ induced by an one--to--one transformation 
$h(\cdot)$; see for instance \cite{Zacks:1971}[Theorem 5.1.1] or \cite{Linkov:2005}[Theorem 5.2.1]. 
(The restriction to one--to--one transformations is not necessary.) 
The question is whether  $\vartheta:=h(\theta)$ 
is the original parameter of another distribution within the same distributional class and how can we get to this distribution by a variable transformation $g(X)$ on the random variable $X$, so that the MLE $\widehat\vartheta_n$ of $\vartheta$ can be computed, not by employing the likelihood of the new distribution but by direct reference to the MLE $h(\widehat \theta_n)$. If this is so we will see that certain methodological implications emerge that considerably simplify associated goodness--of--fit (GOF) procedures.     

There exist already a couple of  well known examples.  Let us start with scale families of distributions, 
and in this regard assume that $X$ has a distribution function (DF) $F(x):=F_c(x)$ such that $F_c(x)=F_0(x/c)$, for some $c>0$, and for some fixed DF $F_0$. 
In this case the family--preserving variable transformation is given by $Y=g(X)$ with $g(X)=a X$, and the new parameter corresponding to the transformed random variable $Y$ is $\sigma=h(c)$ with $h(c)=a c$, for each $a>0$. Then it may be shown that the MLE  satisfies $\widehat\sigma_n(Y_1,...,Y_n)=a \widehat c_n(X_1,...,X_n)$, where $Y_j=a X_j, \ (j=1,...,n)$; see e.g. Problem 3.1 (b), p. 211, in \cite{Lehmann:Casella:1998}. 
The other well known example emerges in location--scale families whereby the DF of $X$ is given by $F_{\theta}(x)=F_0((x-\delta)/c))$, for some $\theta=(\delta,c)$, 
($\delta\in \mathbb R$ and $c>0$), with family--preserving transformation $g(X)=a X+b$, ($a>0$ and $b\in \mathbb R$). Then the new parameter is $\vartheta=(\mu,\sigma)$, where
$\vartheta=h(\delta,c)$ with $h(\delta,c)=(a \delta+b,ac)$, and the MLE satisfies  
$\widehat\vartheta_n(Y_1,...,Y_n)=(a \widehat\delta_{n}+b,a\widehat c_n)$, where $Y_j=a X_j+b, \ (j=1,...,n)$.     
Such estimators that mimic the behavior of corresponding parameters under specific data--transformations are labeled as equivariant estimators; see Definition 2.5 in Chapter 3 of \cite{Lehmann:Casella:1998} for a rigorous definition. 

Here we go beyond these clearly linear cases and pinpoint reasonably rich families of distributions within which, specific distributions differ in more essential ways than mere location and/or scale. Specifically assume that the distribution of $X>0$, depends on a scale parameter $c>0$ and a shape parameter $\kappa>0$, and that the corresponding DF is such that \be \label{fam} F_{c,\kappa}(x)=F_0\left(\left(\frac{x}{c}\right)^\kappa\right),\ee for some fixed DF $F_0(\cdot)$. Specific families of distributions arise by different ``kernel'' DF choices $F_0(\cdot)$ and include:
\begin{enumerate}
    \item The Weibull distribution with kernel DF, $F_0(x)=1-{\rm{e}}^{-x}, \ x>0$.
    \item The Pareto type I distribution with kernel DF, $F_0(x)=1-x^{-1},\ x>1$.
    \item The Frechet distribution with kernel DF, $F_0(x)={\rm{e}}^{-1/x}, \ x>0$.
    \end{enumerate}

In this paper we identify a root--type transformation as a family--preserving variable transformation $g(\cdot)$ within families satisfying \eqref{fam}, as well as the corresponding parameter transformation  $h(\cdot)$ and show that  the MLE satisfies  
\be \label{eqvar} \widehat\vartheta_n=h(\widehat\theta_n),\ee 
for the specific choice of this transformation. 
This is done in Section \ref{sec_2}. Then in Section \ref{sec_3} we investigate the impact that \eqref{eqvar} has on certain goodness--of--fit tests (GOF) for the families of distributions that satisfy \eqref{fam}. At the same time we also provide a non--technical and selective review of the related literature, including  cases of test invariance for linear family--preserving variable transformations. The article concludes with a small simulation study and 
real--data examples in the case of GOF testing for the Weibull distribution in Section \ref{sec_4}.

\section{Distribution--preserving data transformations and equivariant estimators}\label{sec_2} In view of \eqref{fam}, notice that 
\be \label{fam1} X \sim F_{c,\kappa} \implies aX^{1/b}\sim F_{ac^{1/b},b \kappa},\ee  
for each $a,b>0$. Hence the one--to--one root--type transformation $X\mapsto aX^{1/b}$ is a family--preserving variable transformation for each specific family $F_0$ satisfying \eqref{fam}, and maps $(c,\kappa)$ to $(ac^{1/b},b \kappa)$. The following proposition shows that \eqref{eqvar} holds for the MLE in such families of distributions, with variable transformation $g(X)=a X^{1/b}$ and parameter transformation  $h(c,\kappa)=(ac^{1/b},b \kappa)$.  

\begin{prop} \label{prop}
\label{thm_2_5} Let $(X_j, \ j=1,...,n)$ be independent copies of $X\sim F_{c,\kappa}$, 
and assume that the {\rm{DF}} of $X$ satisfies \eqref{fam}. Assume further that the density 
$f_{c,\kappa}$ corresponding to $F_{c,\kappa}$ exists, and also that the {\rm{MLE}} 
$(\widehat c_n,\widehat \kappa_n)$ of the parameter  $(c,\kappa)$ exists. 
Then the {\rm{MLE}} corresponding to $a X^{1/b}\sim F_{ac^{1/b},b \kappa}$ with the sample 
$(aX^{1/b}_j, \ j=1,...,n)$ is given by $(a \widehat c^{1/b}_n ,b\widehat \kappa_n)$. 
\end{prop}

\noindent {\sc Proof.} 
We first calculate from \eqref{fam} the density corresponding to $F_{c,\kappa}$ as    
\be \label{den}
f_{c,\kappa}(x)=\frac{{\rm{d}}F_{c,\kappa}(x)}{{\rm{d}}x}
=\frac{{\rm{d}}F_0((\frac{x}{c})^\kappa)}{{\rm{d}}x}=\frac{\kappa}{c}
 \left(\frac{x}{c}\right)^{\kappa-1}f_0\left(\left(\frac{x}{c}\right)^\kappa\right),
\ee
where $f_0(x):=f_{1,1}(x)=\frac{{\rm{d}}F_0(x)}{{\rm{d}}x}$. 

Then by straightforward calculations, it follows that the likelihood function corresponding to $X$ 
with the sample $(X_j, \ j=1,\ldots,n)$ 
is given by  
\be \label{mle}
{\cal{L}}(X_1,\ldots,X_n;c,\kappa)=\frac{\kappa^n}{\prod_{j=1}^n X_j} \prod_{j=1}^n \left(\frac{X_j}{c}\right)^\kappa  
f_0\left(\left(\frac{X_j}{c}\right)^\kappa\right). 
\ee
Recall {\color{blue}{now}} that $Y:=aX^{1/b}\sim F_{ac^{1/b},b\kappa}$, and denote by $(\widehat d_n, \widehat \eta_n)$ the MLE for $Y$
with the sample $(Y_j=aX^{1/b}_j, \ j=1,\ldots,n)$
, i.e. 
\[
 \max_{ac^{1/b},b\kappa} {\cal{L}}(Y_1,\ldots,Y_n ; ac^{1/b}, b\kappa) = {\cal{L}}(Y_1,\ldots,Y_n ;\widehat d_n,\widehat \eta_n). 
\]
Then, by the definition of $(\widehat c_n,\widehat \kappa_n)$ and $(\widehat d_n,\widehat \eta_n)$, we observe that for any $(c,\kappa)$, the following hold true:
\begin{align*}
& \prod_{j=1}^n \left(\frac{b X^{1-1/b}_j}{a} \right) \cdot {\cal{L}}(X_1,\ldots,X_n;c,\kappa) \\
&=  \frac{(b\kappa)^n}{{\prod_{j=1}^n Y_j}} \prod_{j=1}^n \left(
\frac{Y_j}{a c^{1/b}}
\right)^{b\kappa} f_0\left(\left(\frac{Y_j}{a c^{1/b}}\right)^{b\kappa} \right) \\
&= {\cal{L}}(Y_1,\ldots,Y_n ; ac^{1/b},b\kappa) \\
&\le {\cal{L}}(Y_1,\ldots,Y_n; \widehat d_n , \widehat \eta_n) \\
&= \prod_{j=1}^n \left(\frac{b X^{1-1/b}_j}{a} \right)
\frac{(\widehat \eta_n/b)^n }{\prod_{j=1}^n X_j} 
\prod_{j=1}^n \left( \frac{X_j}{(\widehat d_n /a)^b}\right)^{\widehat \eta_n/b} 
f_0 \left( \left( \frac{X_j}{ (\widehat d_n / a)^b} \right)^{\widehat \eta_n / b} \right) \\
&= \prod_{j=1}^n \left(\frac{b X^{1-1/b}_j}{a} \right) \cdot {\cal{L}}(X_1,\ldots,X_n ; (\widehat d_n /a)^b,\widehat \eta_n/b )\\
&\le \prod_{j=1}^n \left(\frac{b X^{1-1/b}_j}{a} \right) \cdot {\cal{L}}(X_1,\ldots,X_n ; \widehat c_n,\widehat \kappa_n ).
\end{align*}
Now putting $(c,\kappa)=(\widehat c_n,\widehat \kappa_n)$ in the first line of the argument above we obtain  
\[
 {\cal{L}}(X_1,\ldots,X_n; \widehat c_n, \widehat \kappa_n) = {\cal{L}}(X_1,\ldots,X_n; (\widehat d_n /a)^b,\widehat \eta_n/b ).
\]
Hence, $(\widehat d_n,\widehat \eta_n)=(a \widehat c_n^{1/b},b\widehat \kappa_n)$ should be concluded. \bewend

It should be pointed out that there exist other data--transformations that lead to  parameter--free tests, such as the log--linear  transformation $\kappa(\log X-\log c)$ that turns a Weibull variate to a variate following a standard extreme--value distribution. The difference with the root--type transformation suggested herein is that this transformation is not distribution--specific, but instead it applies to rich families of distributions.  In this connection, more general distributions with extra parameters may be included in our framework of families of distributions. For instance  the Burr type XII distribution belongs to the class of distributions defined by \eqref{fam}, with DF $F_0(x)=1-(1+x)^{-\xi}, \ \xi> 0$, and thus also satisfies Prop. \ref{prop}. Specifically if $X$ follows a Burr type XII distribution with parameters $(c,\kappa,\xi)$ then $aX^{1/b}$ follows a  Burr type XII distribution with parameters $(ac^{1/b},b\kappa,\xi)$. Two other well known classes are the exponentiated Weibull distribution with $F_0(x)=(1-{\rm{e}}^{-x})^\xi, \ \xi>0$, and the generalized gamma distribution with density given by  \eqref{den} where $f_0(x)=(\Gamma(\xi))^{-1}x^{\xi-1}{\rm{e}}^{-x}, \ \xi>0$. 

To the best of our knowledge the most general class of distributions that satisfy \eqref{fam} and Prop. \ref{prop} is the ``interpolating family'' 
of distributions on $(0,\infty)$ recently introduced by \cite{Sinner:etal:2023}. As it will be seen below Prop. \ref{prop}  allow us to  carry out a much simpler GOF test that refers to a subfamily of the family under test whereby two (out of two, three or even four) parameters have been removed.

\section{Parameter--free test procedures}\label{sec_3}
It should be noted that the notion of equivariant estimators, and invariant tests (to be discussed below) is a recurring theme in Statistics, with Chapter 3 (resp. Chapter 6) of the classical treatment of \cite{Lehmann:Casella:1998} (resp. \cite{Lehmann:Romano:2005}) devoted to such estimators (resp. tests). 
Nevertheless the approach in \cite{Lehmann:Casella:1998} is mostly based on estimation optimality, which is not always relevant in the context of 
GOF testing. By way of example, smooth tests of fit originally  introduced by \cite{Neyman:1937}, and more recently studied by \cite{Ledwina:1994} 
and in the monograph by \cite{Rayner:etal:2009}, are intrinsically related to moment estimators which are more often than not less efficient than 
other estimators. In fact \cite{Klar:2000} points out that the method of moments is the only meaningful estimation method in the context of smooth tests of fit. 
Another point in case is made by \cite{drost:etal:1990} who argue that rather than estimation optimality, robustness or more precisely luck of it, is important in the context of GOF testing.
Moreover, the test optimality approach often adopted by \cite{Lehmann:Romano:2005}[Chapter 6] 
(see also \cite{Vexler:hutson:2023}) is not even feasible in any reasonably wide context of testing, such as GOF testing with unspecified parameters. This is noted in \cite{Lehmann:Romano:2005}[\S14.1],  and is formally stated and shown by \cite{Jassen:2000}  and \cite{Escanciano:2009}.
(It should be pointed out however that within the narrow context of testing a distribution against a specific alternative, likelihood ratio tests applied on maximal invariants lead to optimal invariant tests). 
Therefore, it appears that the methodological implications of equivariant estimators on GOF testing with estimated parameters have not been put forward beyond the simple linear transformation case of scale or location--scale families, and even in those cases they have not been sufficiently emphasized. 

On the basis of the preceding discussion we motivate our current parameter--free procedures by starting again with simple scale families of distributions with DF $F(x)=F_0(x/c)$. As already mentioned in the Introduction, the MLE of the scale parameter $c$ satisfies $\widehat c_n(a X_1,...,aX_n)=a \widehat c_n(X_1,...,X_n)$. As a result,  any  GOF test for such families of distributions, say $T_n(X_1,...,X_n)$, that depends on  $X_j$ only through $\widehat Y_j=X_j/\widehat c_n$, $( j=1,...,n)$, is scale invariant, i.e. it satisfies 
\be \label{scaleinv} T_n(aX_1,...,aX_n)=T_n(X_1,...,X_n),\ee for each $a>0$, and  consequently the test $T_n$ may be applied by assuming without loss of generality that we are testing for $F=F_0$ with $c=1$. (If $X>0$, this is a special case of Prop. \ref{prop} for $\kappa=1$, and eqn.  \eqref{scaleinv} follows from eqn. \eqref{shapeinv} below for $b=1$). The by far most popular such testing problem is that of testing for exponentiality, and the reader is referred to the review articles of 
\cite{Henze:Meintanis:2005} and \cite{allison:etal:2017}, for scale invariant tests for exponentiality.   

One level up are location--scale families whereby the corresponding location--scale invariant test results by considering the MLE--standardized observations $\widehat Y_j=(X_j-\widehat \delta_n)/\widehat c_n, \ (j=1,...,n)$, and analogously satisfies 
\be \label{locinv} T_n(aX_1+b,...,aX_n+b)=T_n(X_1,...,X_n),\ee for each $a>0$ and $b\in\mathbb R$, implying that $T_n$ can be performed by setting $(\delta,c)=(0,1)$. The monograph of \cite{Thode:2002} focuses on the most popular testing problem here, i.e. that of testing for normality, but without particular reference to test invariance. 
(Nevertheless most normality tests are routinely applied on $(X_j-\widehat \delta_n)/\widehat c_n, \ (j=1,...,n)$, with $\widehat\delta_n$ (resp. $\widehat c_n$) being the sample mean (resp. sample standard deviation), which automatically implies the location--scale invariance stated in \eqref{locinv}). 
On the other hand, \cite{epps:2005} considers GOF tests for general location--scale families with implicit reference to invariance. It should be mentioned that the methodology in \cite{epps:2005} is confined to tests utilizing the empirical characteristic function as their main tool, which might seem as a somewhat less well known approach, but with minor modifications the location--scale invariance argued in Epps (2005) applies more generally to any given GOF test. Moreover, the paper itself as well as earlier (see e.g., \cite{epps:pulley:1983}, \cite{epps:1993}) 
and  subsequent works (see e.g., \cite{Hall:etal:2013}),  including some of the papers to be referenced herein  (see e.g., \cite{Meintanis:Swanepoel:2007}, \cite{Meinitanis:etal:2015}) make a good case about using  the empirical characteristic function  for GOF testing, instead of more standard tools such as  the empirical DF. 

We will briefly digress from univariate distributions, to discuss the very important case of multivariate data. In this connection we note that location--scale equivariance has been extended to vectorial observations as ``affine--equivariance''. The most relevant context for affine equivariant estimators and affine invariant tests is that of (multivariate) elliptical distributions  and an excellent discussion of such estimators and tests may be found in \cite{Hallin:Jureckova:2012}. 
In fact \cite{Hallin:Jureckova:2012} argue  that affine--invariant tests, i.e. tests that for arbitrary dimension $p\geq 1$ satisfy, $T_n(AX_1+b,...,AX_n+b)=T_n(X_1,...,X_n)$, for any non--singular $p \times p$ matrix $A$ and any vector $b\in\mathbb R^p$, should be based on the Mahalanobis distances between the observed $p$--dimensional vectors $(X_j, \ j=1,...,n)$, a point also made by \cite{Henze:2002} in the context of testing for multivariate normality, and by \cite{Meinitanis:etal:2015} for the more general elliptically symmetric stable--Paretian distribution.

Let us return now to our main problem. In this regard,   Prop. \ref{prop} entails that the MLE   estimator $(\widehat c_n,\widehat \kappa_n)$ of $(c,\kappa)$ mimics the equivariance properties of the respective parameters  shown in \eqref{fam1}, i.e. 
\be \label{c} \widehat c_n(a X^{1/b}_1,...,a X^{1/b}_n)=a\:\widehat c^{1/b}_n(X_1,...,X_n),\ee  
and 
\be \label{kappa}  \widehat \kappa_n(a X^{1/b}_1,...,a X^{1/b}_n)=b\:\widehat \kappa_n(X_1,...,X_n),\ee  
for each $a,b>0$.  As a result,  any GOF procedure that depends  on $(X_j, \ j=1,...,n)$ only via 
\be \label{Y} \widehat Y_j=\left(\frac{X_j}{\widehat c_n}\right)^{\widehat\kappa_n}, \ j=1,...,n, \ee
satisfies 
\be \label{shapeinv} T_n(aX^{1/b}_1,...,aX^{1/b}_n)=T_n(X_1,...,X_n),\ee for each $a,b>0$, 
and consequently and without loss of generality, we may perform the test by assuming $c=\kappa=1$. Clearly, in view of the distributional invariance of $aX^{1/b}$ figuring in \eqref{fam1}, test invariance as illustrated by \eqref{shapeinv} feels like a desirable, even natural, property within the families of distributions satisfying \eqref{fam}.  On the practical level eqn. \eqref{shapeinv} implies that a potentially much simpler test may be invoked for testing families satisfying \eqref{fam}, such as in the case of the Weibull distribution where any test for exponentiality applied on $(\widehat Y_j, \ j=1,...,n)$ can be used.           

\begin{rem} Notice that the data--transformation figuring in eqn. \eqref{Y} is the inverse of the distribution--preserving root--type transformation $g(X)$ shown in \eqref{fam1}, and thus not--surprisingly, its application has a stabilizing effect on the estimators. Specifically, by replacing  in \eqref{c}--\eqref{kappa}, $(a,b)$ by $(\widehat c_n^{-\widehat \kappa_n}, \widehat \kappa^{-1}_n)$ we see easily that   
\[ \widehat c_n(\widehat Y_1,...,\widehat Y_n)=\widehat c_n\left(\left(\frac{X_1}{\widehat c_n}\right)^{\widehat\kappa_n},...,\left(\frac{X_n}{\widehat c_n}\right)^{\widehat\kappa_n}\right)
=1,\]  
and 
\[  \widehat \kappa_n(\widehat Y_1,...,\widehat Y_n)=\widehat \kappa_n\left(\left(\frac{X_1}{\widehat c_n}\right)^{\widehat\kappa_n},...,\left(\frac{X_n}{\widehat c_n}\right)^{\widehat\kappa_n}\right)=1.
\]
\end{rem}

\begin{rem} The invariance properties figuring in \eqref{scaleinv} and \eqref{locinv} are not restricted to the MLE alone. Other estimators, such as moment estimators or estimators based on order statistics, may also result in test statistics that satisfy these properties provided that the estimators under discussion  satisfy the equivariance properties referred to in the Introduction. A particular case of equivariant estimators with minimum risk are the Pitman estimators of location and scale; see \cite{Zacks:1971}[\S7.2], \cite{Lehmann:Casella:1998}[\S3.1] and \cite{Linkov:2005}[\S3.1-3.2].   
\end{rem}


In this connection, and before closing this section we wish to emphasize that the invariance properties figuring in \eqref{scaleinv}, \eqref{locinv}, and \eqref{shapeinv}, do not imply that our test procedures reduce to the case of simple hypotheses with corresponding parameters known. In fact parameter estimation generally does have an effect on the distributional properties of the tests, and test invariance only means that these distributional properties do not involve the actual true  values of the unknown parameters being estimated; see for instance \cite{Meintanis:Swanepoel:2007}.            
In the next section we illustrate the performance of some GOF tests for exponentiality that are employed in order to test for the Weibull distribution with both parameters unknown.


\section{Monte Carlo and real--data}\label{sec_4}
In this section we study the finite--sample performance of a few tests for the Weibull distribution with DF,  $F_{c,\kappa}(x)=1-\exp\{-(x/c)^\kappa\}$, and unknown parameter $(c,\kappa)$. Recall that if the tests are applied on $(\widehat Y_j, \ j=1,\ldots,n)$ as defined in \eqref{Y} with $(\widehat c_n,\widehat \kappa_n)$ being the MLE, then we can set $c=\kappa=1$, and consequently we may invoke any test for unit exponentiality (see item 1 in Section \ref{sec_1}). From the plethora of available tests we consider the Anderson--Darling test based on 
\[
 {\rm{AD}}_{n} =-n -\frac{1}{n} \sum_{j=1}^n 
(2j-1) (\log{Z_{(j)}}+\log (1-Z_{(n+1-j)})), 
\]
where $Z_{(j)}=1-{\rm{e}}^{-\widehat Y_{(j)}} \ (j=1,...,n)$ with  $\widehat Y_{(1)}\leq \widehat Y_{(2)}{\leq} \ldots \leq \widehat Y_{(n)}$ being the ordered statistics. The AD$_n$ test is often the most powerful test among the classical tests based on the empirical DF. We also include the test 
of \cite{Henze:Meintanis:2002} based on 
 \[
 {\rm{HM}}_{n} = \frac{1}{n} \sum_{j,k=1}^n 
\frac{1+(\widehat Y_j+\widehat Y_k+2)^2}{(\widehat Y_j+\widehat Y_k+1)^3}-2 \sum_{j=1}^n 
\frac{\widehat Y_j+2}{(\widehat Y_j+1)^2} +n, 
\]
which is amongst the best performing exponentiality tests in the comparison studies of \cite{Henze:Meintanis:2005} and \cite{allison:etal:2017}, 
and a smooth test of fit for the exponential distribution (see \cite{Rayner:etal:2009}[\S6.3] given by 
\[
{\rm{RB}}_{n}=\frac{1}{n}  \left(\sum_{j=1}^n L_2(\widehat Y_j)\right)^2+\frac{1}{n}\left(\sum_{j=1}^n L_3(\widehat Y_j)\right)^2,
\]
where 
\[ 
L_2(z)=1-2z +z^2/2, \ L_3(z)=1-3z+3 z^2/2-z^3/6, \] 
are the Laguerre polynomials of orders 2 and 3.

We consider tests of size $\alpha$, and for a given sample size $n$ we calculate the test statistics based on a large number $M$ of Monte Carlo samples and obtain the quantile  corresponding to $1-\alpha$. Specifically for each sample of size $n$ and each test statistic, say $T$, we generate observations from a Weibull distribution with a fixed combination of $(c,\kappa)$, then we calculate the MLE $(\widehat c_n, \widehat \kappa_n)$ and  obtain the value of the test statistic $T_m$ based on the transformed sample $\widehat Y_j= (X_j/\widehat c_n)^{\widehat \kappa_n}, \ j=1,...,n$. By iterating this procedure for $m=1,2,...,M$, we obtain the critical value of the test statistic as the $1-\alpha$ quantile of the empirical distribution of   $(T_m, \ m=1,...,M)$. By the test invariance articulated in the previous section we only need to draw samples from the Weibull distribution with $(c,\kappa)=(1,1)$, i.e. by sampling from the unit exponential distribution. Nevertheless we examined the three test statistics by sampling from a Weibull distribution with several combinations of $(c,\kappa)$ and indeed our conclusion for a parameter--free test statistic  was confirmed as the resulting critical values remained stable regardless of the actual value of $(c,\kappa)$ employed. 

The actual Monte Carlo was performed with sample size $n=50, 100, 150$ and $n=200$, with $M=100,000$ iterations at significance level  $\alpha=0.1,0.05$ and $\alpha=0.01$, and the resulting critical values are reported in Table \ref{table:criticalpt}. The figures in Table   \ref{table:criticalpt} show that convergence to the asymptotic distribution is faster for the AD$_n$ and HM$_n$ tests, while the smooth test RB$_n$ is somewhat slower to reach its limit distribution.     

\begin{table}[hbtb]
\begin{center}
\caption{Critical values for AD$_{n}$, HM$_{n}$, and RB$_{n}$, at significance level $\alpha$}  
\begin{tabular}{|l|ccc|ccc|ccc|}\hline \label{table:criticalpt}
       & \multicolumn{3}{|c|}{AD$_{n}$} & \multicolumn{3}{|c|}{HM$_{n}$} & \multicolumn{3}{|c|}{RB$_{n}$}\\ \hline
       $n\setminus\alpha$  & $0.1$ & $0.05$ & $ 0.01$   & $0.1$ & $0.05$ & $0.01$ & $0.1$ & $0.05$ & $0.01$\\  \hline
 $n=50$ & $0.629$ & $0.750$ & $1.027$ & $0.036$ & $0.047$ & $0.076$ & $0.922$ & $1.374$ & $4.558$ \\
  $n=100$  & $0.623$ & $0.755$ & $1.027$ & $0.037$ & $0.048$  & $0.078$  &$1.123$& $1.797$  & $7.456$\\ 
 $n=150$ & $0.634$ & $0.755$ & $1.041$ & $0.037$ & $0.049$ & $0.079$& $1.269$& $2.092$ & $9.249$\\
   $n=200$  & $0.634$ & $0.757$ & $1.050$ & $0.038$ & $0.049$  & $0.078$  & $1.383$ & $2.313$ & $9.914$\\ \hline 
  \end{tabular}
\end{center}
\end{table}

\begin{figure}[hbpt]
\begin{center}
\includegraphics[width=0.8\textwidth]{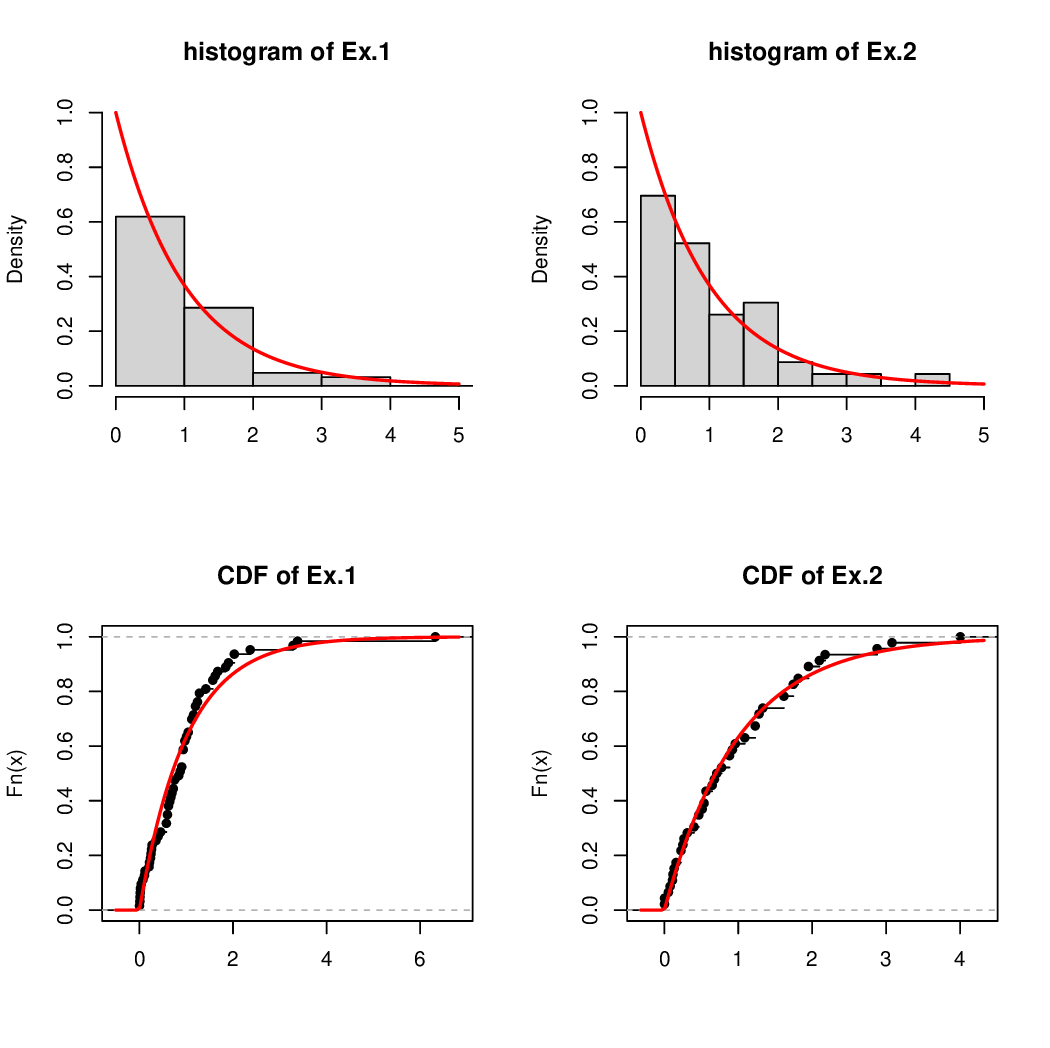}
\end{center}
\caption{Histogram and distribution function of the standardized strength data 
of glass fiber (Ex.1 and Ex.2). The standard exponential density and distribution function are 
superimposed on the corresponding graphs (solid curved lines).  
}\label{fig:1}
\end{figure}

Using this methodology we apply the three tests on two real--data sets of sizes $n=63$ (Ex.1) and $n=46$ (Ex.2)
employed by \cite{Smith:Naylor:1987} (see Table 1, p.359 for the data). 
The data correspond to experimental measurements on the strength of glass fiber of length $1.5$ cm (Ex.1) and $15$ cm (Ex.2). 
\begin{table}
\begin{center}
\caption{Actual values of test statistics (``act.''), and critical values at significance level 0.1, 0.05, 0.01}
\begin{tabular}{|l|cccc|cccc|cccc|}\hline \label{table:actual}
       & \multicolumn{4}{|c|}{AD$_{n}$} & \multicolumn{4}{|c|}{HM$_{n}$} & \multicolumn{4}{|c|}{RB$_{n}$}\\ \hline
        & act. & $0.1$ & $0.05$  & $0.01$ & act. & $0.1$ & $0.05$ & $0.01$& act. & $0.1$ & $0.05$ & $0.01$\\  \hline
 Ex.1 & $\underline{1.241}$ & $0.629$ & $0.747$ & $1.031$& $\underline{0.104}$ & $0.036$ & $0.047$ & $0.076$& $\underline{1.458}$ & $0.981$ & $1.496$ & $5.820$ \\
 Ex.2 & $\underline{0.329}$ & $0.631$ & $0.752$ & $1.024$ & $\underline{0.021}$ & $0.036$ & $0.046$  & $0.077$ &$\underline{0.602}$ &$0.904$& $1.349$  & $4.277$\\  \hline 
  \end{tabular}
\end{center}
\end{table}
The values of the MLE are $(\widehat c_n, \widehat \kappa_n)=(5.781,1.628)$ for Ex.1, and 
$(\widehat c_n, \widehat \kappa_n)=(5.147,1.230)$ for Ex.2. 
From these estimated values, we obtain $(\widehat Y_j, \ j=1,...,n)$, and calculate each of the three test statistics, whose actual values 
are reported in Table \ref{table:actual} (underlined figures). 
Before further analysis we also considered the Kolmogorov-Smirnov (KS) test with the data $\widehat Y_j$, and thereby obtained p-values $0.1078$ for Ex.1 and $0.9473$ for Ex.2. These are consistent with the values reported in Table 4 (Ex.1) and Table 5 (Ex.2), of the real--data analysis in  \S5 of \cite{Wu:etal:2021}. 
Thus at significance level $10\%$, and on the basis of the KS test we can not reject the hypothesis of exponentiality of the standardized data $(\widehat Y_j, \ j=1,...,n)$, which in turn can be interpreted to imply that the source data  $(X_j, \ j=1,...,n)$ might have originated from a Weibull distribution with the MLE estimates as parameters. On the other hand,  
the exponentiality of the data of Ex.1 is rejected by the AD$_n$ and the HM$_n$ tests, and only the smooth test at significance levels $5\%$ and $1\%$ finds no evidence to reject this hypothesis. At the same time, the exponentiality for the data of Ex.2 is supported by all three tests uniformly over the values of $\alpha$ considered. 
The corresponding histogram and distribution function plots that are shown in Figure \ref{fig:1} further corroborate our results. 
Therefore, there is strong evidence in favour of a Weibull distribution for the data of Ex.2, while the corresponding conclusion for the data of Ex.1 should be questioned and occasional non--rejection may be due to low power. For instance, the KS test is often the least powerful amongst the classical tests based on the empirical DF.


{}

\end{document}